\documentclass[12pt,a4paper]{amsart} %parametri reqno siirtää equation numberssit oikealle
\usepackage[T1]{fontenc}
\usepackage[utf8]{inputenc}
\usepackage{amsmath, amssymb, amsthm, latexsym}
\usepackage[margin=3cm]{geometry}

\theoremstyle{plain}
\newtheorem{thm}{Theorem}
\newtheorem{lemma}[thm]{Lemma}

\newtheorem{proposition}[thm]{Proposition}
\newtheorem{fact}[thm]{Fact}

\theoremstyle{definition}
\newtheorem{definition}[thm]{Definition}
\newtheorem{remark}[thm]{Remark}

\newtheorem{notation}[thm]{Notation}
\newtheorem{question}[thm]{Question}

\numberwithin{thm}{section}
\numberwithin{equation}{section}

\renewcommand{\a}{\alpha}
\renewcommand{\b}{\beta}
\newcommand{\e}{\varepsilon}

\renewcommand{\o}{\omega}

\newcommand{\R}{\mathbb{R}}
\newcommand{\Q}{\mathbb{Q}}
\newcommand{\C}{\mathbb{C}}
\newcommand{\abs}[1]{\lvert#1\rvert}
\newcommand{\norm}[1]{\lVert#1\rVert}
\newcommand{\raj}{\restriction}
\newcommand{\ket}[1]{|#1\rangle}
\newcommand{\dom}{\operatorname{dom}}
\newcommand{\rng}{\operatorname{rng}}
\newcommand{\nat}{\mathbb{N}}
\newcommand{\Z}{\mathbb{Z}}

\begin{document}
\title[Eigenvectors, Approximations and the Feynman Propagator]{On Eigenvectors, Approximations and the Feynman Propagator}

\date{\today}
\author{{\AA}sa Hirvonen}
\address{Department of Mathematics and Statistics\\
University of Helsinki\\
P.O.Box 68\\
00014 University of Helsinki\\
Finland}
\email{asa.hirvonen@helsinki.fi}
\author{Tapani Hyttinen}
\address{Department of Mathematics and Statistics\\
University of Helsinki\\
P.O.Box 68\\
00014 University of Helsinki\\
Finland}
\email{tapani.hyttinen@helsinki.fi}

\begin{abstract}
Trying to interpret B. Zilber's project on model theory of quantum mechanics we
study a way of building limit models from finite-dimensional approximations. 
Our point of view is that of metric model theory, and we develop a method of taking ultraproducts of unbounded operators. We first calculate the Feynman propagator for the free particle 
as defined by physicists as an inner product $\langle
x_{0}\vert K^{t}\vert x_{1}\rangle $ of the eigenvector $\vert
x_{0}\rangle $ of the position operator with eigenvalue $x_{0}$ and
$K^{t}(\vert x_{1}\rangle )$, where $K^{t}$ is the time evolution operator.
However, due to a discretising effect, the eigenvector method does not work as expected, and without heavy case-by-case scaling, it gives the wrong value. We look at this phenomenon, and then complement this by showing how to instead calculate the kernel of the time evolution operator (for both the free particle and the harmonic oscillator) in the limit model.
We believe that
our method of calculating these is new.
\end{abstract}

\maketitle

\section{Introduction}

In this paper we study an approximation technique for quantum mechanics. We develop a technique for approximating the $L_2(\R)$ model by finite-dimensional spaces, and illustrate the use of it by calculating the kernels of the time-evolution operator for the free particle and harmonic oscillator.
There is a tradition in physics of using such finite-dimensional approximations, but it is not always clear in what sense these approximations happen. The core of our paper is the approximation theory which justifies how and when this can be done, and in particular, points out a problem in a first, direct approach. Our approximation theory is based on metric ultraproducts, and we develop a way to treat (nice enough) unbounded operators in these ultraproducts, so that the technique can be used. The problem we discover is a discretising effect tied to the use of Gauss sums, which causes direct calculations with eigenvectors to give the wrong answer. We then show how to remedy the situation, by instead of the propagator considering the kernel of the time evolution operator.

Our interest in the subject stems from B. Zilber's talk 
on the use of model theory in quantum physics, which he gave at
Jouko V\"a\"an\"anen's 60th birthday meeting.
We got interested and tried to read
\cite{Zi} in which Zilber calculates the Feynman propagator
for the free particle and the harmonic oscillator.
But we failed to follow the calculations. However, the paper
suggests various ways of constructing models that resemble
the Hilbert space operator model used in quantum physics
for a single particle.
In addition,
we noticed that by using number theory, in the case of
the free particle, the Feynman propagator is easy to calculate in
many of these
models (the case of the harmonic oscillator is harder).
Our first straight forward calculations gave incorrect values for the
propagator, leading us to study various ways of scaling
the operators
and renormalising the results
to get the right answer. This seems to be common in physics, see e.g.\ \cite{TC}.
Looking closer at these
scalings we noted that the use of eigenvectors (Dirac delta-functions)
causes a clumping effect on the propagator. This eventually lead us to calculating
the kernel instead of the propagator. We look at this phenomenon in section \ref{sec:calc}.

The approach used by Zilber is based on finite dimensional approximations
of the standard model from quantum mechanic, and Zilber is not the only
one that has considered these. The idea of finite dimensional approximations goes back to Weyl \cite{We} and Schwinger
\cite{Sc} (in various approaches,
there are differences in details such as what  the dimensions
one considers are and how the negative part of the space in which the particle
lives is treated, but, although  technically important,
we overlook these difference in this introduction).
The first question in these considerations is: Do the approximations really
approximate the operators we  are interested in and if they do,
in what sense do they approximate them? Neither Weyl nor Schwinger defines precisely what is meant by the approximation, and the more modern expositions seem to rely either on approximating bounded operators or operators with a pure point spectrum.

In \cite{Go} and \cite{AGK} the approximation theory is built by nonstandard analysis (and thus ultraproducts), in \cite{DVV} the approximation is in standard mathematics. Both \cite{AGK} and \cite{DVV} approximate the 
Hamiltonian of the harmonic oscillator based on the well-understood
spectral theory of this operator which tells pretty much everything
about the operator. E.g., in \cite{AGK} the approximation theory
is developed for quasicompact operators  and although
the Hamiltonian $H$ of the harmonic oscillator is not
quasicompact,
by using the spectral theory one can see that
e.g., $(H-I)^{-1}$ is, and thus the theory can be used.
However,  for operators like the position operator or
the Hamiltonian of the free particle these methods do not  work.

We develop a theory of  approximations  that
applies to all operators that arise in quantum mechanics (i.e., also operators with continuous spectrum).
The point here is that we want a theory  that has a chance to be used
to study  also operators of which we  do not in advance have a good
understanding (or whose behaviour is known to be complicated).
We demonstrate these
possibilities by calculating
the kernel of the time evolution operator for both the free particle and
the harmonic oscillator and do this without using any spectral information of the Hamiltonian.
And this brings us to the main theme of this paper:
the difference between the Feynman propagator and the kernel
of the time evolution operator in their behaviour in the approximations.

In \cite{DHV} the authors give the coefficients of the
Hamiltonian of the harmonic oscillator in the approximations but
they do not tell how they  calculate them.
However, in
\cite{TC}, J. Tolar and G. Chadzitaskos calculate the Feynman propagator
for the free particle much the same way as we do
using Siegel's formula.
The main differences between
what we do and what is done in \cite{TC}  are the following:
First of all, in \cite{TC} the authors do not seem to have any theory of
approximations. This may also be the reason why they
do the calculations only for the free particle, while
using the theory of approximations, we can deal with the harmonic oscillator
case as well.
And finally
we will show that
without taking averages, this eigenvector approach does not really work.

The main motivation for studying finite dimensional approximations is that
the theory of operators in finite dimensional Hilbert spaces  is much
simpler than the theory in infinite dimensional spaces.
E.g., if $P_{N}$ is a 
decently well behaved %(esim. nilpotenteille operaattoreille t??m?? ei p??de)
operator in a finite dimensional space
$H_{N}$, then $H_{N}$ has an orthonormal basis that consists
of eigenvectors of $P_{N}$. One does not in general have the same luxury
in the infinite dimensional spaces
(in the case of the Hamiltonian of the harmonic
oscillator one does have this luxury). E.g., the eigenvector
$\vert x_{0}\rangle $ for position $x_0$ does not exist
(in $L_{2}(\R )$, that is, in our metric ultraproduct it exists).
Thus calculations in finite dimensional spaces  may be easier
and by  \L os's theorem for (our variant of) metric ultraproducts,
these determine what happens
in the metric ultraproduct.

In quantum mechanics physical observables are described by self adjoint
operators on a complex Hilbert space, which usually
is of infinite dimension. The unit
sphere of the space 
corresponds to the possible states of the system and the possible
outcomes of measurements are the eigenvalues of the
operators. If a self adjoint operator has non-degenerate eigenvalues,
the corresponding eigenvectors will be orthogonal with real
eigenvalues. 
Physicists work
under the assumption that there not only are eigenvalues
with unique eigenvectors, but that the
eigenvectors of any of the operators considered span the whole
space. This is called 'inserting a
complete set of states'.
So they work in a Hilbert space spanned by the eigenvectors
corresponding to the possible outcomes of measurements. Linear
combinations of these are considered 'superpositions', states in which
the observable is not yet determined until it is measured (upon which
the state 'collapses' onto one of the eigenstates). The coefficients
in the linear combinations determine the probability that the state
will collapse onto the corresponding eigenstate when the observable is
measured.

The evolution of the system in time is described by
unitary time evolution operators $K^{t}$, i.e., if the state at time
$0$ is $\phi$, then at time $t$ it is  $K^{t}(\phi )$.
The time evolution operator is found as a solution
to the Schr\"odinger equation $i\hbar\frac{\partial}{\partial t}K^t=HK^t$,
where $H$ is the so called Hamiltonian, i.e., the operator for energy.

In the case of a single free particle in space,
the assumption of a complete set of states clashes with the convention
of working in a separable Hilbert space.
In fact, it is easy to see that, e.g., the standard position operator (see below)
does not have any eigenvectors.
One approach to overcome this problem is to
work in a rigged Hilbert space (i.e., Gelfand triple) but
mathematical physicists tend to
override this problem by simply forgetting the
eigenvectors. Instead they study regions of the spectrum within
which the measurements can fall. The motivation for this is that one
in reality cannot measure point values anyway, but only, e.g., whether a
particle is inside a given area (the detector) of space. Thus 
notions involving eigenvectors are replaced with methods based on the
spectral theorem (for unbounded self adjoint operators), 
and instead of considering specific positions in space, mathematical
physicists just use wave functions giving the position of the particle
a probability distribution over the possible places. Let us look closer at how
this affects the calculations of the 
Feynman propagator for the free particle.

The way physicists interpret the {\it propagator} $K(x,y,t)$ (for a time
independent Hamiltonian) is that
it expresses the probability amplitude of the
particle  travelling form point $x$ to point $y$ in time $t$. If we
have eigenstates $|x\rangle$, $|y\rangle$ corresponding to the
positions, this probability amplitude will be given by the inner
product of $|y\rangle$ and $K^t|x\rangle$ where $K^t$ is the time
evolution operator, i.e.,
$$K(x,y,t)=\langle y|K^t|x\rangle.$$
This is the basic idea behind methods used by physicists. It is also
the starting point for both Zilber's and our calculations.

If one does not assume eigenvectors, the propagator is defined as the
kernel of the integral representation of the time evolution operator
$$(K^t\psi)(y)=\int_{\R} K(x,y,t)\psi(x)dx.$$
This approach is described, e.g., by Lukkarinen \cite{Lu}. Note that mathematically
the propagator and the kernel 
are different things but they are used the same way when
one wants to calculate the movement of the particle and thus
they should have the same value (otherwise the two approaches
give different physics).

In both approaches the first thing to do is to find
the time evolution operator $K^{t}$.
In the approach described in \cite{Zi}, spectral theory is used
but Lukkarinen relies more on Fourier analysis. 
Below in the spaces $H_{N}$ 
we find the time evolution operator using the spectral theory method,
only now the needed spectral theory is trivial
since the dimensions of the spaces are finite.

Then in the  approach used in \cite{Zi}, an integral value for the propagator
is calculated using the Feynman path integral. Lukkarinen, on the
other hand, shows by a straight forward
calculation that
for Schwartz test functions $f$,
$$K^{t}(f)(x)= \int_{\bf R} dk\ e^{i2\pi k(x-y)-it(1/2)(2\pi k)^2} $$
(here Lukkarinen has done some scaling and the units are chosen
so that both $\hbar$ and the mass of the particle are one).
Finally, in \cite{Zi} a value for the integral is found by a leap of faith
but Lukkarinen argues as follows:
Let us denote by $K^{t}_{*}$ the operator defined using
the correct answer as the kernel. So one needs to
show that $K^{t}=K^{t}_{*}$.
For this, for all positive reals $\e$, one defines
operators
$$K^{t}_{\e}(f)(x)= \int_{\bf R} dk\ e^{i2\pi k(x-y)-(it+\e)(1/2)(2\pi k)^2} $$
(this is called regularising $K^{t}$).
By using dominant convergence theorem, Fubini
and the known values for Gauss integrals, for
Schwartz test functions $f$ it is possible to show that
$$K^{t}(f)=\lim_{\e \rightarrow 0}K^{t}_{\e}(f)=K^{t}_{*}(f).$$
This suffices.

To calculate the propagator using finite-dimensional approximations,
there are a few things one needs to consider. The operators considered
are the position operator $Q$ and the momentum operator $P$. In the
standard space $L_2(\R)$ these are $Qf(x)=xf(x)$ and $Pf(x)=-i\hbar (df/dx)(x)$
(here $\hbar =h/2\pi$, where $h$ is the Planck constant
and we restrict to the case where the space in which the particle
lives is of dimension one). Their commutator is $[Q,P]=i\hbar$. 
However, in no finite dimensional Hilbert space can one find
self adjoint operators $Q$ and $P$ satisfying the commutation
relation.
Also, the operators $Q$ and $P$ are not bounded, so defining metric
ultraproducts of them poses another problem.
Thus as a starting point for our model we use the exponentials
$U^{t}=e^{itQ}$ and $V^{t}=e^{itP}$ of the operators. These will satisfy
the Weyl commutator law $V^{w}U^{t}=e^{i\hbar tw}U^{t}V^{w}$
for many $t$ and $w$.

The Hilbert space model for $P$ and $Q$ for a single particle in one
dimensional space  
(from the point of view of, e.g.,
the Feynman propagator the dimension assumption is w.l.o.g.)
is governed by two theorems.
The first one is Stone's theorem which states that
if for all $t\in\R$, $U^{t}$ is a unitary operator
in a complex Hilbert space $H$, $t\mapsto U^{t}$ is continuous
(i.e. $t\mapsto U^{t}(x)$ is continuous for all $x\in H$) and
$U^{t+t'}=U^{t}U^{t'}$, then there is a unique self adjoint
operator $Q$ such that $U^{t}=e^{itQ}$ for all $t\in\R$
(typically the self adjoint operator is not bounded
and so by $e^{itQ}$ we mean the continuation to the whole space).
We will call the maps $t\mapsto U^{t}$ as above
continuous unitary representations of $(\R ,+)$.

The other central theorem is the Stone-von Neumann theorem. It states
that the class of Hilbert spaces $H$ with
continuous unitary representations $t\mapsto U^{t}$
and $t\mapsto V^{t}$ of $(\R ,+)$ with the property
$V^{w}U^{t}=e^{i\hbar tw}U^{t}V^{w}$ is essentially categorical
(in the sense of model theory i.e. 'unitarily')
in every cardinality.
However, we want to point out
that if $t\mapsto U^{t}$
and $t\mapsto V^{t}$ are as above and $a$ is, e.g., a positive real
then so are $t\mapsto U_{*}^{t}$
and $t\mapsto V_{*}^{t}$ where $U_{*}^{t}=U^{at}$
and $V_{*}^{t}=V^{a^{-1}t}$ and that this kind of scaling may have
an effect in
calculations.

Now since neither our finite dimensional approximations
nor our full limit model 
are the standard
model, it is thus clear that we will have to give up either the unitary
representations of $(\R,+)$, or the Weyl commutator law has to
fail. What happens is that in our finite dimensional models we have
unitary representations, but the commutation relation only holds
partially. Thus  we have operators $Q'$ and
$P'$ whose exponentials $U^t$ and $V^t$  are, but they are not the
real position and momentum operator, and we will look at how they
relate to their correct counterparts in our calculations in section \ref{sec:calc}.
Further, in moving to the metric ultraproduct we lose continuity of
the representations.
Also, since the Weyl commutator law does not
hold in full in the finite dimensional models, in the metric ultraproduct
it holds
only in a dense set of $t$s and $w$s. 

However, we will further develop
our techniques of taking metric ultraproducts to find a smaller limit
model for the finite dimensional approximations inside the
ultraproduct and in this submodel 
the  Weyl commutator law holds and the continuous
unitary representations remain continuous, thus producing a copy of
the standard model that can be approximated by finite dimensional
spaces in the sense of \L os's theorem (some extra requirements are needed).
To find this copy, 
we introduce an isometric embedding $F$ of
$L_{2}(\R )$ into our limit model. We
show that all operators relevant to the physics are mapped
naturally to a 
metric ultraproduct of operators in 
the finite dimensional approximations.
This allows us 
to decide not only properties of the limit model
but also properties of the model used  by  physicists
by calculating in the finite dimensional approximations.

Our calculation of the propagator in section \ref{propagator} gives a result that is either scaled or 0, depending on divisibility issues. In an earlier version of this paper we claimed that these don't matter, as we 'almost always' are in a divisible case, and thus only need to take care of the scaling along with a 'renormalisation' factor. However, this is not true, as in the limit model any given position will have continuum many corresponding eigenvectors, and getting divisibility depends on the choice of representative. This is the reason one needs to look at kernels instead. In section \ref{scalings} we demonstrate why scaling cannot resolve the problem with eigenvectors (Dirac deltas). 

However, scaling also opens up new possibilities for applications. In section \ref{bohr} we show how using a radically different scaling, our model gives Bohr's model for the hydrogen atom.

\section{Using metric model theory}

In this section we show how to approximate the $L_2(\R)$ model of quantum mechanics with finite dimensional spaces. This includes developing methodology for treating unbounded operators in ultraproducts.

Our starting point will be a metric version of \L os's theorem stating that a formula is true in an ultraproduct model if and only if it is true in almost all (in the filter sense) index models. We will embed the $L_2(\R)$ model of quantum mechanics into a metric ultraproduct of finite dimensional Hilbert spaces, thus making precise what we mean by approximations: The finite dimensional spaces approximate the $L_2(\R)$ model in the sense that the properties inherited by $L_2(\R)$ from the ultraproduct (the properties preserved by the embedding) are precisely the ones true in almost all the approximating finite-dimensional spaces.

Ultraproducts of Banach spaces were originally defined in \cite{BDCK} and gained new attention via the rise of continuous logic (see, e.g., \cite{mtcts}). The basic idea is to start with an indexed collection of Banach spaces and an ultrafilter on the index set. For elements of this product, one then looks at the ultralimits of their coordinatewise norms, and restrict attention to only the elements, for which this is finite. Finally, one mods out the infinitesimals.

For uniformly bounded families of bounded operators, the ultraproduct is defined the classical way, i.e., by just taking the equivalence class of the sequence of coordinatewise values. However, we wish to study also operators, that do not have a uniform bound in our index models, so we need to develop a methodology for this first. The result will be a partially defined unbounded operator.

\subsection{Ultraproducts and unbounded operators}\label{ultraproducts}

We study a method for taking ultraproducts of operators that are not uniformly bounded in the index models. In our applications the index models will be finite dimensional Hilbert spaces (where, of course, no unbounded operators exist), and to ease notation we assume the operators we study are bounded, just not uniformly so. However, this assumption is, in fact, without generality, as our other assumptions will cut out subspaces where the operator is bounded.

Before stating the theorem, let us recall the definition of ultraproduct of Hilbert spaces and fix some notation:

Given Hilbert spaces $H_i$, $i\in I$ and an ultrafilter $D$ on $I$, we can define the metric ultraproduct $H$ of the spaces $H_i$:

  First let $H_*$ be the set of those $f\in\prod_{i\in I}H_i$ for which there is some $M\in\o$ such that $\{i\in I:\norm{f(i)}<M\}\in D$. Then let $H=H_*/\sim_D$, where $\sim_D$ is the equivalence relation
  $$
  f\sim_D g \textrm{ if for all }M<\o,\, \{i\in I:\norm{f(i)-g(i)}<\frac{1}{M}\}\in D.
  $$
  For $f\in\prod_{i\in I}$, we write $f/D$ for the equivalence class of $f$ under this equivalence relation.

Technically our models are many sorted, as we have one sort for the
complex numbers, and sorts for the Hilbert space. If one wants to work
in the framework of continuous logic, one needs the Hilbert spaces to
be bounded, and then one divides them up into $N$-radius balls, $N<\o$, all of which are added as sorts to the structures
(this can be avoided by using positive bounded formulae instead, see
\cite{HI} for details).

In order to simplify the notation, if for all $i\in I$,
$\xi_{i}\in H_{i}$ is given, then by $(\xi_{i})/D$
we mean $f/D$, where $f\in\prod_{i\in I}H_{i}$ is such that
for all $i\in I$, $f(i)=\xi_{i}$. When convenient, we denote
$f/D$ also by $\lim_{D}\xi_{i}$. And if $\xi_{i}\in\C$
for all $i\in I$, by $\lim_{D}\xi_{i}$ we mean analogously
the complex number $r$ such that for all $\e >0$,
$$\{i\in I : \ \abs{\xi_{i}-r} <\e\}\in D,$$
if such an $r$ exists.

It is straightforward to show that for uniformly continuous functions with a common modulus of uniform continuity, the usual definition of function in an ultraproduct gives a well defined function with the same modulus of uniform continuity. The main novelty of this section is that in some cases one can construct ultraproducts of functions without a common modulus of uniform continuity:

\begin{thm}\label{K-konstruktio}
  Let, for each $i\in I$, $H_i$ be a complex Hilbert space and $P_i$ a bounded 
  operator on $H_i$ (where the bound may vary with $i$). Further assume there are complete subspaces $H_i^k$ (possibly $\{0\}$), for all $k<\o$, such that
  \begin{enumerate}
  \item if $k\neq l$, then $H_i^k$ and $H_i^l$ are orthogonal to each other,
  \item $P_i(H_i^k)\subseteq H_i^k$,
  \item for all $k<\o$, there is $0<M_k<\o$ such that for all $i\in I$ and $x\in H_i^k$
    $$
    \frac{1}{M_k}\norm{x}\leq\norm{P_i(x)}\leq M_k\norm{x}.
    $$
  \end{enumerate}
  Then if $D$ is an ultrafilter on $I$, there is a closed subspace $K$ of the metric $D$-ultraproduct of the spaces $H_i$ where we can define the ultraproduct of the operators $P_i$ as an unbounded operator $P$ satisfying
  \begin{enumerate}
  \item on a dense subset of $K$, $P(f/D)=(P_i(f(i)))_{i\in I}/D$ and
  \item if for $n<\o$, $f_n/D\in\dom(P)$ and both $(f_n/D)_{n<\o}$ and $(P(f_n/D))_{n<\o}$ are Cauchy sequences, and $(f_n/D)_{n<\o}$ converges to $f/D$, then $P$ is defined at $f/D$ and $P(f/D)=\lim_{n\to\infty}P(f_n/D)$.
  \end{enumerate}
\end{thm}

\begin{proof}
  Let $H_i^k$ and $P_i$ be as in the theorem, and let, for each $n<\o$, $K_i^n$ be the (complete) subspace of $H_i$ generated by $\bigcup_{k\leq n}H_i^k$.

  Consider the metric ultraproduct $H$ of the spaces $H_i$. Here,  to guarantee a proper treatment of our
partially defined operators, we need to add sorts for the subspaces
$K_i^n$ (which further can be split up into balls as described above). 

For all $n<\o$,
we say that $f\in\prod_{i\in I}H_{i}$ is \emph{$n$-guarded} if there is
$M<\o$ such that
$$\{i\in I:\ f(i)\in K_{i}^{n},\ \norm{f(i)}\leq M\}\in D.$$
We say that $f$ is \emph{guarded} (or \emph{$\o$-guarded})
if it is $n$-guarded for some $n<\o$.
We let $K_{n}$ be the complete subspace of $H$ generated by
$\{ f/D: f\textrm{ is $n$-guarded}\}$.
We write just $K$ for $K_{\o}$.
Notice that if $f$ is ($n$-)guarded, then so is $(P_{i}(f(i)))_{i\in I}$. 

Next, for each $f/D\in K$, we say that $(f_{n}/D)_{n<\o}$ is the canonical
Cauchy sequence of $f/D$ if for all $n<\o$, $f_{n}/D$
is the orthogonal projection of $f/D$ to $K_{n}$.
It is easy to see that $(f_{n}/D)_{n<\o}$
is a Cauchy sequence and that
it converges to $f/D$.
In fact one can choose the functions $f_{n}$ so
that for all $i\in I$, $f_{n}(i)$ is the  projection of $f(i)$
to $K^{n}_{i}$.
Notice that if $(f_{n}/D)_{n<\o}$ and $(g_{n}/D)_{n<\o}$
are canonical Cauchy sequences for $f$ and $g$, respectively, then
$((\a f_{n}+\b g_{n})/D)_{n<\o}$ is the canonical Cauchy sequence
for $\a f+\b g$ for all complex numbers $\a$ and $\b$.

Now we can define
$P$ in $K$ as follows (it is easy to see that $P$
is well-defined; we do not define $P$  outside $K$).
If $f$ is guarded, we let
$P(f/D)=(P_i(f(i)))_{i\in I}/D$
and otherwise we let
$(f_{n}/D)_{n<\o}$ be the canonical
Cauchy sequence of $f/D$ and let
$P(f/D)=\lim_{n\rightarrow\infty}P(f_{n})/D$
assuming $\lim_{n\rightarrow\infty}\norm{P(f_{n})/D}<\infty$
(the limit exists)
and if $\lim_{n\rightarrow\infty}\norm{P(f_{i})/D} =\infty$,
$P(f/D)$ is undefined (thus $P$ may be unbounded
in $K$). It is easy to see that $P$ is an unbounded operator in $K$ with the desired properties.
\end{proof}

\begin{remark}
If $f$ is not guarded, it is essential to use the
canonical Cauchy sequence to define $P(f/D)$.
To see why, let $I=\o$ and for each $N\in\o$, let $H_N$ be a finite dimensional Hilbert space spanned by the orthogonal basis $\{u(0),\dots, u(N)\}$. Let $H_N^k=\langle u(k)\rangle$ and $P(u(n))=nu(n)$. Now consider $x=(\sum_{k=0}^N\frac{1}{k^2}u(k))_{N<\o}/D$. If we denote $\xi_{N,n}=\sum_{k=0}^{\min\{n,N\}}\frac{1}{k^2}u(k)\in H_N$, then $((\xi_{N,n})_{N<\o}/D)_{n<\o}$ is the canonical Cauchy sequence for $x$. But if we instead look at the sequence built from elements $\xi'_{N,n}=\xi_{N,n}+\frac{1}{n}u(n)$, then the $f'_n=(\xi'_{N,n})_{N<\o}/D$ are $n$-guarded and converge to $x$. However, in each $H_N$ with $N>n$ we have $P(\xi'_{N,n})=P(\xi_{N,n})+u(n)$, so although $f_n$ and $f'_n$ converge to the same element, $P(f_n)$ converges to $P(x)$, and $\{\norm{P(f'_n)}:n<\o\}$ is bounded, the sequence $P(f'_n)$  is not a Cauchy sequence.
\end{remark}

Next we look at the exponential $e^{iP}$ of $iP$
under the assumption that $P_{j}$
is self adjoint for all $j\in I$.

For all  $j\in I$, we define
$e^{iP_{j}}$ the usual way:
For $x\in H_{j}$,
$$e^{iP_{j}}(x)=x+\sum_{n=1}^{\infty}(iP_{j})^{n}(x)/n!.$$
Then $e^{iP_{j}}$ is a unitary operator  in $H_{j}$
and  for all $m<\o$,  $e^{iP_{j}}(K_{j}^{m})\subseteq K_{j}^{m}$.
Then for any ($n$-)guarded $f=(f_{j})_{j\in I}$,
$e^{iP}(f)=(e^{iP_{j}}(f_{j}))_{j\in I}$ is ($n$-)guarded
and so we can define $e^{iP}(f/D)=(e^{iP_j}(f(j)))_{j\in I}/D$.
Finally, for an arbitrary $f/D\in K$, we
let $(f_{m}/D)_{m<\o}$ be the canonical
Cauchy sequence of $f/D$ and define
$e^{iP}(f/D)=\lim_{m\rightarrow\infty}e^{iP}(f_{m}/D)$
(here we could have chosen any Cauchy sequence of guarded
functions in place of the canonical one and got the same
result). The following fact is an immediate consequence
of the standard theory of metric ultraproducts
applied to the spaces $K^{n}_{j}$:

\begin{fact}\label{fact1.3}
\begin{enumerate}
\item  $e^{iP}$ is the metric ultraproduct of the operators
$e^{iP_{j}}$ restricted to $K$.
\item\label{toinen} For all guarded $f$,
$$e^{iP}(f/D)=f/D+\sum_{n=1}^{\infty}(iP)^{n}(f/D)/n!.$$
\end{enumerate}
\end{fact}

In general, we now have a method to apply \L os's theorem also for
operators that are unbounded in $H$, as long as we consider only parameters from
$K$. The basic vocabulary for this consists of the addition in the Hilbert space,
the inner product from $(H_{j})^{2}$ to $\C$,
the scalar multiplication $\C\times H_{j}\rightarrow H_{j}$,
the field structure on $\C$, the norm $|\cdot|$ in $\C$ (or the
complex conjugation), 
the complex exponentiation and constants for $i$, $h$ and $\pi$. We
will also have unitary operators $U^t$ and $V^t$ defined on the whole
space and the 'true' operators from quantum mechanics defined on (the
sort/sorts corresponding to) $K$.

\subsection{Building a model for quantum mechanics}\label{building}

The starting point of our model (roughly following Zilber) is a space
with unitary operators $U^t$ and $V^t$ satisfying the Weyl commutator relation 
$$V^{s}U^{t}=e^{ist\hbar}U^{t}V^{s}$$ 
(Zilber uses
$V^{s}U^{t}=e^{2\pi ist\hbar}U^{t}V^{s}$).
These will be defined on all of the
space. They will agree with the quantum mechanical operators $e^{itQ}$
and $e^{itP}$ on $K$ for a dense set of $t$s.

As the commutator relation leaves
room for scaling, we pick two 
positive reals $a'$ and $b'$ for the scaling. We will look closer at
this scaling along the way - first at the requirements we get on $a'$
and $b'$ to make our model fit the requirements of quantum mechanics, and
later on the effects the scaling choice has on calculations.

For all natural numbers $N\in\o$, let $H_{N}$
be a vector space over the complex numbers $\C $ with
basis $\{ u(x)\vert\ x\in\nat ,
\ x<N\}$ ($H_0=\{0\}$, and it is included just for convenience).

Now, technically, for $N<N'$ the basis for $H_{N}$ is a subset of that of $H_{N'}$, but this  is not a natural way of looking at these spaces, we do it this way in order to simplify the notation.

We make Hilbert spaces out of the vector spaces $H_{N}$ by defining
an inner product $\langle\ \cdot\ |\ \cdot\ \rangle $ so that
$\langle u(x)|u(x')\rangle =0$ if $x\ne x'$ and otherwise
$\langle u(x)|u(x')\rangle =1$ (and for complex numbers $c$ and $d$,
$\langle cu(x)|du(x)\rangle =\overline{c}d$).

When we take roots of complex numbers $c=e^{i2\pi r}$, $0\le r<1$,
we mean the principal ones, i.e., $c^{m/n}=e^{i2\pi mr/n}$.

By $q=q_{N}$ we denote the complex number
$e^{i2\pi /N}$ and define
$$v(x)=(1/N)^{1/2}\sum_{y=0}^{N-1}q^{xy}u(y) ,$$
the 'Fourier transform' of $u(x)$.
Notice that
for $x\ne x'$,
$$\langle v(x)|v(x')\rangle =(1/N)\sum_{y=0}^{N-1}e^{i2\pi zy/N}$$
for some integer $z\ne 0$ with $|z|<N$ and
thus $\langle v(x)|v(x')\rangle =0$. Similarly
one can see that $\langle v(x)|v(x)\rangle =1$, and so
$\{ v(x)\vert\ x\in\nat ,
\ x<N\}$ is another orthonormal basis of $H_{N}$.
Furthermore, similarly, one can see that
$$u(x)=(1/N)^{1/2}\sum_{y=0}^{N-1}q^{-xy}v(y),$$
since $$u(x)\mapsto (1/N)^{1/2}\sum_{y=0}^{N-1}q^{-xy}v(y)$$
and $$v(x)\mapsto (1/N)^{1/2}\sum_{y=0}^{N-1}q^{xy}u(y)$$
are inverses of each other.

For all real numbers $t$, let $t_0^{u}=a't/2\pi$ and
$t_0^{v}=b't/2\pi$. 
Also for all real numbers $t$, we define operators $U^{t}$ and $V^{t}$
as follows: For all $0\le x<N$,
$U^{t}(u(x))=q^{xt_0^{u}}u(x)$ and
$V^{t}(v(x))=q^{xt_0^{v}}v(x)$. Notice that if
$t_0^{v}$ is a natural number, 
then
$$V^{t}(u(x))=u(x-t_0^{v})$$
where the sum is taken 'modulo $N$'
i.e. $-y=N-y$ if $N\ge y>0$. This is because

$$u(x-t_0^{v})=(1/N)^{1/2}\sum_{y=0}^{N-1}q^{-(x-t_0^{v})y}v(y)$$
and
$$V^{t}(u(x))=V^{t}((1/N)^{1/2}\sum_{y=0}^{N-1}q^{-xy}v(y))=
(1/N)^{1/2}\sum_{y=0}^{N-1}q^{-xy}q^{t_0^{v}y}v(y).$$

So now
$U^{t}$ and $V^{w}$ are unitary onto operators,
$u(x)$ is an eigenvector of $U^{t}$ with eigenvalue
$e^{ixa't/N}$ and $v(x)$ is an eigenvector of $V^{w}$
with eigenvalue $e^{ixb'w/N}$
and assuming that both $t_0^{u}$ and $w_0^{v}$ are natural numbers,
$V^{w}U^{t}(u(x))=q^{xt_0^{u}}u(x-w_0^{v})$
and $U^{t}V^{w}(u(x))=q^{(x-w_0^{v})t_0^{u}}u(x-w_0^{v})$
i.e. $V^{w}U^{t}=q^{w_0^{v}t_0^{u}}U^{t}V^{w}=e^{ia'b'tw/2\pi N}U^{t}V^{w}$.

Now it is time to look at the Weyl commutator law.
We want to have $$V^{w}U^{t}=e^{i\hbar tw}U^{t}V^{w},$$
so we need to require that $a'b'=Nh$.

For the 'infinite' scaling $N$, there are only a few choices
that  give any chances of getting the right results for the propagator,
and out of those only one seems to give reasonable physics
(in section \ref{bohr} we look at another scaling giving a circular universe). Thus
we use the natural even distribution of $N$ between the
two operators,  i.e., we let $a'=\sqrt{N}a$ and $b'=\sqrt{N}b$
where $a$ and $b$ are positive real numbers such that
$ab=h$.

It will turn out, that the right scaling to use in calculations is $a=1$ and $b=h$, but we will give our basic definitions with respect to general choices of $a$ and $b$ to be able to look closer at the reason.

\begin{notation}\label{notaatio}
Having fixed the 'infinite' part of the scaling, we define $t^{u}=at/2\pi$ and $t^{v}=bt/2\pi$. So $t^u=t_0^u/\sqrt{N}$, $t^v=t_0^v/\sqrt{N}$, and thus
$$
U^{t}(u(x))=q^{x\sqrt{N}t^{u}}u(x)
$$
and the $U^t$-eigenvalue of $u(x)$ is $e^{ixat/\sqrt{N}}$. Accordingly
$$
V^{t}(v(x))=q^{x\sqrt{N}t^{v}}v(x)
$$
and the $V^t$-eigenvalue of $v(x)$ is $e^{ixbt/\sqrt{N}}$. Also
$$
V^t(u(x))=u(x-\sqrt{N}t^v)
$$
when $\sqrt{N}t^v$ is an integer. Further note that if we define operators $Q'_N$
and $P'_N$ in $H_N$ by $Q'_Nu(x)=xaN^{-1/2}u(x)$ and
$P'_Nv(x)=xbN^{-1/2}v(x)$, then $U^t=e^{itQ'_N}$ and $V^t=e^{itP'_N}$,
but $Q'_N$ and $P'_N$ do not behave like the standard position and
momentum operators (we cannot have a constant commutator in a finite dimensional space). We can also define the operators $U_2^t=e^{it(Q'_N)^2}$
and $V_2^t=e^{it(P'_N)^2}$ with eigenvalues $e^{it(ax)^2/N}$ and
$e^{it(bx)^2/N}$, respectively, and these will be important in our calculations.
\end{notation}

In addition to the choice of 'infinite' scaling and choosing the values for $a$ and $b$, there is a third scaling, that occurs when we start calculating the propagators.
To use Gauss summation formulae, we will need some
parameters to have rational values. But this can easily be achieved by
suitably choosing units for time, length or mass.

Now we take a metric ultraproduct of these structures.
Let $D$ be an ultrafilter on $\o$
such that for all $n>0$, $X_{n}\in D$, where
the set $X_{n}$ consists of those $N\in\o$ for which $\sqrt{N}$ is
a natural number and
divisible by $n$. We say that some claim is true in $H_{N}$
for almost  all $N$ if the set of those $N$ for which the claim is true
belongs to  $D$. And, in fact, often we do not mention this 'for almost all'
at all. Also in definitions, it is enough that the definition works
for almost all $N$.

We thus have the following:

\begin{proposition}
The model $H$ is a complex Hilbert space with unitary operators $U^t$ and $V^t$ such that for $ab=h$ and all $t$ and $w$ such that $t^u=at/2\pi$ and $w^v=bt/2\pi$ are rational numbers (and thus  $\sqrt{N}t^u$ and
$\sqrt{N}w^v$ are integers for almost all $N$), 
$$V^{w}U^{t}=e^{i\hbar tw}U^{t}V^{w}.$$

Furthermore, the model has eigenvectors:
assuming $a^{-1}x$ is a rational number,
we have (e.g.) that $u(x)=(u(a^{-1}\sqrt{N}x)\vert\ N\in X)/D$ is
an eigenvector of each $U^{t}$
with eigenvalue $e^{itx}$.
Similarly,
(e.g.) $v(x)=(v(\sqrt{N}x)\vert\ N\in X)/D$ is an eigenvector
of each $V^{t}$
with eigenvalue $e^{itbx}$.
\end{proposition}

However, it is worth noting,
that the functions $(t,x)\mapsto V^{t}(x)$
and $(t,x)\mapsto U^{t}(x)$ do not have a
modulus of uniform continuity (the unitary operators $U^t$ and $V^t$
themselves, of course, have one). So in the full ultraproduct the operators $U^t$ and
$V^t$ do not form continuous unitary representations of $(\R,+)$. Of
course we know this, since otherwise the Stone-von Neumann would
guarantee we are still essentially in the standard $L_2(\R)$ model. 

By looking at embeddings of the $L_2(\R)$ model we can, however, find a submodel of our ultraproduct, where we do have unitary representations, as well as the 'true' quantum mechanics operators.

\subsection{Embeddings of the $L_2(\R)$ model}\label{embedding}

In this section we show how to embed the standard $L_2(\R)$ model into our ultraproduct model $H=\prod_{N\in\o}H_N/D$ preserving all essential operators.

There are (at
least) two reasons for looking at such  
embeddings. One is that we want our methods to say something about the
physics predicted by the standard model. Thus our model should offer a
representation of the standard model, and what we do is add something
outside this. So in a sense we add eigenvectors to the standard model
and get something (a lot) bigger. Note, however, that what we do is
not the 'insertion of a complete set of states' that the physics
literature talks of. We do add eigenvectors for all possible
positions (and momenta), but these do not span all of the model.

The other reason we need an embedding, is 
that we need the standard model to say something about our
calculations. In section~\ref{oscillator} we will use a factorisation of the time-evolution operator to calculate
the kernel of the harmonic oscillator. In the finite-dimensional
spaces $H_N$ the product we calculate will not give the right operator, 
but since the factorisation holds in the standard model, the
product we use will give the right operator in a piece of our model
containing (the image of) the standard model.

We start with a set of scaled embeddings to see how the scaled
operators $U^t$ and $V^t$ behave.

\begin{definition}
For fixed scaling parameters $a$ and $b$ (with $ab=h$), 
let 
$F=F^{a}_{b}$ be the following isometric embedding of $L_{2}$ into $H$:
In a dense set of ('nice') functions $f\in L_{2}$ (e.g., $C^\infty_c$,
the set of compactly supported smooth functions)
let  $F(f)=(F_{N}(f)\vert\ N<\o)/D$, where for $N>1$
\begin{align*}
F_{N}(f)  =  & \sum_{n=0}^{(N/2)-1}a^{1/2}N^{-1/4}f(naN^{-1/2})u(n)+\\
 & \sum_{n=N/2}^{N-1}a^{1/2}N^{-1/4}f((n-N)aN^{-1/2})u(n).
\end{align*}
As $F$ is isometric, it can be extended to all of $L_2(\R)$.
\end{definition}

So $F_N$ looks at functions restricted to the interval 
$[-a\sqrt{N}/2,a\sqrt{N}/2]$ and maps them into $H_N$ by 
cutting the interval into pieces of length $a/\sqrt{N}$ and 
looking at the values of $f$ at the starting points of these intervals.
This determines $F$ everywhere and is isometric.

It is straightforward to show that these embeddings map the operators $e^{itx}$ and $e^{-it\hbar(d/dx)}$ to the operators $U^t$ and $V^t$ respectively, when $ta/\pi\in\Q$ and $tb/\pi\in\Q$, which is exactly the same criterion on $a,b$ and $t$ as we got from the Weyl commutator relation:

\begin{lemma}
Consider the exponentials of the quantum mechanical operators $Q(f(x))=xf(x)$ and $P(f(x))=-i\hbar\frac{df}{dx}(x)$. When $ab=h $, and $ta/\pi,tb/\pi\in\Q$, then $F^a_b$ maps $e^{itQ}$ and $e^{itP}$ to the operators $U^t$ and $V^t$ respectively (restricted to the image of the embedding).
\end{lemma}

\begin{proof}
For $e^{itQ}$, note that if $f\in C^\infty_c$, then $e^{itQ}(f(x))=e^{itx}f(x)$, so
\begin{align*}
F_N(e^{itQ}f)= & a^{1/2}N^{-1/4}(\sum_{n=0}^{N/2-1}e^{itnaN^{-1/2}}f(naN^{-1/2})u(n)+\\
& \sum_{n=N/2}^{N-1}e^{it(n-N)aN^{-1/2}}f((n-N)aN^{-1/2})u(n))
\end{align*}
and
\begin{align*}
U_N^tF_Nf= &
a^{1/2}N^{-1/4}(\sum_{n=0}^{N/2-1}e^{itnaN^{-1/2}}f(naN^{-1/2})u(n)+\\
 & \hfill\sum_{n=N/2}^{N-1}e^{itnaN^{-1/2}}f((n-N)aN^{-1/2})u(n)).
\end{align*}
Thus $F_N$ maps $e^{itQ}$ to $U_N^t\raj \rng(F_N)$ whenever
$e^{-itaN^{1/2}}=1$, which is when $ta/\pi\in\Q$.

For $e^{itP}$, we have that $e^{itP}f(x)=f(x+t\hbar)$ (see, e.g., Chapter 10 of \cite{Ha}). Thus
\begin{align*}
F_N(e^{itP}f)= & a^{1/2}N^{-1/4}(\sum_{n=0}^{N/2-1}f(naN^{-1/2}+t\hbar)u(n)+\\
& \sum_{N/2}^{N-1}f((n-N)aN^{-1/2}+t\hbar)u(n))
\end{align*}
and if $bt/\pi\in\Q$ (and thus $k:=\sqrt{N}bt/2\pi\in\Z$) (and w.l.o.g. $k<N/2$)
\begin{align*}
V_N^tF_Nf= & a^{1/2}N^{-1/4}(\sum_{n=0}^{N/2-1}f(naN^{-1/2})u(n-\sqrt{N}bt/2\pi)+\\
& \sum_{N/2}^{N-1}f((n-N)aN^{-1/2})u(n-\sqrt{N}bt/2\pi))\\
= &
a^{1/2}N^{-1/4}(\sum_{n=0}^{N/2-k-1}f((n+k)aN^{-1/2})u(n)+\\
& \sum_{n=N/2-k}^{N-k-1}f((n+k-N)aN^{-1/2})u(n)+\\
& \sum_{n=N-k}^{N-1}f((n+k-N)aN^{-1/2})u(n)).
\end{align*}
Now recalling that $ab=h$, we see that $kaN^{-1/2}=t\hbar$, i.e.,
$F_Ne^{itP}$ and $V_N^tF_N$ agree, except when $N/2-k\leq n<N/2$ and by
choosing $N$ large 
enough (compared to the support of $f$) we can make this difference
vanish. So $e^{itP}$ is mapped to $V^t\raj \rng(F)$ when $tb/\pi\in\Q$.
\end{proof}

The proof above shows that the unitary, exponentiated operators are preserved under our embeddings. However, to use our model for calculation, we will need to look at operators of the form $e^{i\alpha(\beta q)^2}$ and $e^{i\alpha(\beta P)^2}$. Thus we will need to show that the real quantum physics operators (whose counterparts will be defined only in $K$), are also preserved. This will require another approach to $F$. 
In order to simplify notations,
we will only consider the scaling
$a=1$. The general case goes exactly the same way.

Now first 
notice that for $C^{\infty}$-functions $f\in L_{2}(\R )$,
$F(f)=\lim_{m\to\infty}F(f_{m})$, where
$f_{m}(x)=f(x)$ if $-m\le x\le m$ and otherwise
$f_{m}(x)=0$ (and $m$ is a natural number).
Thus we can use another definition for $F$:

\begin{definition}
For all rational numbers $q<r$, we let
$\phi_{qr}\in L_{2}(\R )$ be such that
$\phi_{qr}(x)=1/\sqrt{r-q}$ if $q\le x\le r$ and otherwise
$\phi_{qr}(x)=0$. 
For $N<\o$ (such that $\sqrt{N}$ is large and divisible enough),
let $n$ and $m$ be integers such that
$n/\sqrt{N}=q$, $m/\sqrt{N}=r$
and $\vert n\vert ,\vert m\vert <N/2$.
Then let
$F_N(\phi_{qr})=(F_{N}(\phi_{qr}))/D$, where
$F_{N}(\phi_{qr})$ is defined  as it was defined above for
$C^{\infty}_{c}$-functions, i.e.,
if  $q<0$ and $r>0$, let
$$F(\phi_{qr})=(m-n)^{-1/2}(\sum_{k=0}^{m-1}u(k)
+\sum_{k=1}^{-n}u(N-k)),$$
if $q\ge 0$, let
$$F_{N}(\phi_{qr})=(m-n)^{-1/2}\sum_{k=n}^{m-1}u(k)$$
and finally, if $r\le 0$, let
$$F_{N}(\phi_{qr})=(m-n)^{-1/2}
\sum_{k=-m+1}^{-n}u(N-k).$$
By noting, that  $C^{\infty}_{c}$-functions can be approximated
by step functions, it is straightforward to see that this determines an isometric embedding, which is the same as the $F$ defined above (for $a=1$).
\end{definition}

Next we will define another embedding $G$ and then 
show that $G=F$. We will use this to see 
how $F$ maps various operators of $L_{2}(\R )$.

\begin{definition}
Let $\Phi :L_{2}(\R )\rightarrow L_{2}(\R )$
be the isometric automorphism of $L_{2}(\R )$ one gets from extending
the inverse of the 
Fourier transform, i.e., $\Phi$ is the operator such that
for integrable functions $f\in L_{2}(\R )\cap L_1(\R)$,
$$\Phi (f)(x)=\int_{\R}f(y)e^{2\pi i yx}dy.$$ 
For all rational numbers $q<r$, let
$\psi_{qr}=\Phi (\phi_{qr})$.
Then the subspace generated by
the set $S$ of all $\psi_{qr}$ is dense in $L_{2}(\R )$
(since the subspace generated by
the set of all $\phi_{qr}$ is dense in $L_{2}(\R )$)
and thus to define $G$, it is enough to define
$G(\psi_{qr})$ for all  $\psi_{qr}\in S$. We do this as
for $F$ (mapping $\phi_{qr}$) except that we use vectors $v(n)$ in
place of $u(n)$: 
$G(\psi_{qr})=(G_{N}(\psi_{qr}))/D$, where
$G_{N}(\psi_{qr})\in H_{N}$ is defined as follows:
Again let $n$ and $m$ be integers such that
$n/\sqrt{N}=q$, $m/\sqrt{N}=r$
and $\vert n\vert ,\vert m\vert <N/2$.
Then if  $q<0$ and $r>0$
let
$$G_{N}(\psi_{qr})=(m-n)^{-1/2}(\sum_{k=0}^{m-1}v(k)
+\sum_{k=1}^{-n}v(N-k)),$$
if $q\ge 0$, let
$$G_{N}(\psi_{qr})=(m-n)^{-1/2}\sum_{k=n}^{m-1}v(k),$$
and finally, if $r\le 0$, let
$$G_{N}(\psi_{qr})=(m-n)^{-1/2}
\sum_{k=-m+1}^{-n}v(N-k).$$
Notice that if $q<p<r$, then
$G(\psi_{qr})=\a G(\psi_{qp})+\b G(\psi_{pr})$
for suitably chosen real numbers $\a$ and $\b$ (the rescaling coefficients $\a$ and $\b$ are necessary as all the $\psi$-functions have norm 1).
\end{definition}

\begin{lemma}
$G$ is an isometric embedding of $L_{2}(\R )$
into $H$.
\end{lemma}

\begin{proof}
As above,
$\phi_{qr}\mapsto  G(\psi_{qr})$ determines
an isometric embedding $G'$ of $L_{2}(\R )$
into $H$ ($G'$ is similar to $F$, the only difference is that
we have used vectors $v(k)$ in place of $u(k)$).
But then also $G=G'\circ\Phi^{-1}$
is an isometric embedding of $L_{2}(\R )$
into $H$.
\end{proof}

\begin{proposition}\label{prop2.2}
$G=F$. 
\end{proposition}

\begin{proof} It is enough to prove that for all $\psi_{qr}\in S$,
$G(\psi_{qr})=F(\psi_{qr})$.
Notice that the functions $\psi_{qr}$ are 'nice' functions,
and the derivative of $\psi_{qr}$ is $\Phi (2\pi ix\phi_{qr})$.
In order to simplify notations we prove the claim
only for $q=0$ and  $r=1/k$, $k\in\o -\{ 0\}$.

We will show this by looking at the projections $G^m_N$ and $F^m_N$
($m\in\o$) of
$G_N$ and $F_N$ to the subspace $A^m_N$ of $H_N$ generated by 
$\{ u(n)\vert\ 0\le n<m\sqrt{N}\}\cup
\{ u(n)\vert\ N-m\sqrt{N}\le n<N\}$.

For $m\in\o$, let $\psi_{qr}^{m}$ be such that
for all $-m\le x\le m$, $\psi_{qr}^{m}(x)=\psi_{qr}(x)$
and otherwise $\psi_{qr}^{m}(x)=0$.
 Then (for all large enough $N$)
$F^m_N(\psi_{qr})=F_N(\psi^m_{qr})$. Thus, if we denote $F^m=\lim_DF^m_N$,
we have $F^m(\psi_{qr})=F(\psi^m_{qr})$.
Also, as $\psi^m_{qr}\to\psi_{qr}$ and $F$ is isometric, $F(\psi^m_{qr})$
converges to $F(\psi_{qr})$.

Now we need two things to finish the proof: we need the distance
between $G^m_N(\psi_{qr})$ and $F^m_N(\psi_{qr})$ to go to 0 as $\sqrt{N}$
becomes more divisible. And we need
$G^m(\psi_{qr})=\lim_DG^m_N(\psi_{qr})$ to converge to $G(\psi_{qr})$
as $m$ grows.

For the first part, note that 
by presenting $G_{N}(\psi_{qr})$ in the basis of $u(n)$'s,
one gets (denoting $N_r=\sqrt{N}r$):
$$G_{N}(\psi_{qr})=N_{r}^{-1/2}\sum_{k=0}^{N_{r}-1}
N^{-1/2}(\sum_{n=0}^{N-1}e^{2\pi i kn/N}u(n))=$$
$$\sum_{n=0}^{(N/2)-1}N^{-1/4}\sum_{k=0}^{N_{r}-1}(r^{-1/2}
e^{2\pi i(k/\sqrt{N})nN^{-1/2}})N^{-1/2}u(n)+$$
$$\sum_{n=N/2}^{N-1}N^{-1/4}\sum_{k=0}^{N_{r}-1}(r^{-1/2}
e^{2\pi i(k/\sqrt{N})(n-N)N^{-1/2}})N^{-1/2}u(n).$$
Above we used that $N_{r}^{-1/2}=r^{-1/2}N^{-1/4}$
and that $e^{2\pi ik(n-N)/N}=e^{2\pi ikn/N}$.

Further recall that 
$$F^m_N(\psi_{qr})=\sum_{n=0}^{(\sqrt{N}m)-1}
N^{-1/4}(\int_{0}^{r}r^{-1/2}e^{2\pi ixnN^{-1/2}}dx)
u(n)+$$
$$\sum_{n=N-\sqrt{N}m}^{N-1}
N^{-1/4}(\int_{0}^{r}r^{-1/2}e^{2\pi ix(n-N)N^{-1/2}}dx)
u(n).$$

Then note that 
$$\sum_{k=0}^{N_{r}-1}(r^{-1/2}
e^{2\pi i(k/\sqrt{N})nN^{-1/2}})N^{-1/2}$$
converges to
$$\int_{0}^{r}r^{-1/2}e^{2\pi ixnN^{-1/2}}dx$$
uniformly (with error $N^{-1/2}$) for all $n$ such that
$0\le n<m\sqrt{N}$ (the negative side is symmetric).
Combining these observations gives
\begin{equation}\label{eq1}
 \Vert G^m_N(\psi_{qr})-F^m_N(\psi_{qr})\Vert\leq (m\sqrt{N}N^{-1})^{1/2}
\end{equation}
for large enough $N$.

For the second claim, note that as
$G_N(\psi_{qr})-G^m_N(\psi_{qr})$ is orthogonal to $G^m_N(\psi_{qr})$,
also $G(\psi_{qr})-G^m(\psi_{qr})$ is orthogonal to
$G^m(\psi_{qr})$, and thus $\Vert
G(\psi_{qr})-G^m(\psi_{qr})\Vert^2=\Vert G(\psi_{qr})\Vert^2-\Vert
G^m(\psi_{qr})\Vert^2$. Further from \eqref{eq1} it follows that
$G^m(\psi_{qr})=F^m(\psi_{qr})$.
Now as $\psi_{qr}\in L_2(\R)$, for all $\e$ there is
$m_\e$ such that for $m\geq m_\e$,
$\Vert\psi_{qr}-\psi^m_{qr}\Vert<\e$, and $\Vert\psi^m_{qr}\Vert>1-\e$.
Then as both $F$ and $G$ are isometries,
\begin{align*}
\Vert G(\psi_{qr})-G^m(\psi_{qr})\Vert^2
& =\Vert G(\psi_{qr})\Vert^2-\Vert G^m(\psi_{qr})\Vert^2\\
& =\Vert G(\psi_{qr})\Vert^2-\Vert F^m(\psi_{qr})\Vert^2\\
& =1-\Vert \psi_{qr}^m\Vert^2
\end{align*}
which proves that $G^m(\psi_{qr})$ converges to $G(\psi_{qr})$.
And as $F^m(\psi_{qr})$ converges to $F(\psi_{qr})$, this finishes the proof.
\end{proof}

We now turn to the 'true' operators of quantum mechanics:
the position operator $x:f\mapsto xf$
and momentum operator $-i\hbar(d/dx):f\mapsto -i\hbar(df/dx)$
on $L_{2}(\R )$, where $\hbar =h/2\pi$ and $h$ is the Planck
constant. So from now on we will use $x$ both for this operator and also as
a real number. It will be clear from the context which we mean.
Recall that these operators are unbounded and self adjoint.

It can be shown that for Schwartz functions (which are dense
in $L_2$) one has (see, e.g., \cite{Ha}):

\begin{equation}\label{eqstar}
\Phi\circ x\circ\Phi^{-1}=-(i/2\pi )(d/dx),
\end{equation}
which we will use as a definition for $-(i/2\pi )(d/dx)$
(this is a common practice, and then, of course,
$-i\hbar (d/dx)$ means the operator $h(-(i/2\pi )(d/dx))$).

In $H_{N}$ we define an operator $Q_{N}$
so that $Q_{N}(u(n))=nN^{-1/2}u(n)$ if $n<N/2$
and otherwise $Q_{N}(u(n))=(n-N)N^{-1/2}u(n)$
and an operator  $P_{N}$ so that
$P_{N}(v(n))=hnN^{-1/2}v(n)$ again if $n<N/2$
and otherwise $P_{N}(v(n))=h(n-N)N^{-1/2}v(n)$.
To define operators $Q$ and $P$ on a subspace $K$ of $H$
from the operators $Q_{N}$ and $P_{N}$,
we use the method from Section~\ref{ultraproducts}. For this we need to define
the subspaces $H^{k}_{N}$. For $Q$ we let $H^{k}_{N}$
be the subspace generated by
$$\{ u(n)\vert\ kN^{1/2}\le n<(k+1)N^{1/2}\}
\cup\{ u(n)\vert\ kN^{1/2}< N-n\le (k+1)N^{1/2}\}.$$
For $P$ we define the subspaces similarly, but using the $v$-basis, i.e.,
$H^{k}_{N}$ is the subspace generated by
$$\{ v(n)\vert\ kN^{1/2}\le n<(k+1)N^{1/2}\}
\cup\{ v(n)\vert\ kN^{1/2}< N-n\le (k+1)N^{1/2}\}.$$
(So $H^{k}_{N}=\{ 0\}$ if  $k\ge N^{1/2}$.)

Note that this way we get two subspaces $K_Q$ and $K_P$, depending on whether we build $K$ from the $u$- or the $v$-basis of the $H_N$. These are different subspaces of $H$: $K_Q$ has eigenvectors for positions and $K_P$ has eigenvectors for momenta, but it is easy to show that the spaces spanned on one hand by the positional eigenvectors and on the other by the momentum eigenvectors are orthogonal (as the coefficients from the Fourier transform between them become infinitesimal). However, as $F=G$, we know that $L_2(\R)$ embeds into both. How much else fits into their intersection is not known to us:

\begin{question}
What is the intersection $K_P\cap K_Q$?
\end{question}

This is related to the more general question:
\begin{question}
How canonical is $K$ (for a fixed operator)? Can it be essentially altered by the choice of the spaces $H^k_N$ in the construction?
\end{question}

We will study the operator $Q$ a bit later, so for now, when we talk about
$H^{k}_{N}$, $K^{n}_{N}$ or $K_{n}$ we mean those defined for
$P$. Notice that for all $\psi_{qr}\in S$,
$G(\psi_{qr})\in K_{n}$ for every large enough $n$
and thus $\rng(G)\subseteq K$.

\begin{thm}$F\circ -i\hbar (d/dx)\circ F^{-1}=
P\raj F(\dom(-i\hbar (d/dt)))$.
\end{thm}

\begin{proof} By Proposition~\ref{prop2.2}, it is enough to show that for all
$\psi_{qr}\in S$, 
$$G(-i\hbar(d\psi_{qr}/dx))=P(G(\psi_{qr})).$$
As before, we prove this for $q=0$ and $r=1/k$.
We notice that by \eqref{eqstar} above,
$$-i\hbar (d\psi_{qr}/dx)=\Phi (hx\phi_{qr})=
\int_{0}^{r}r^{-1/2}hye^{2\pi iyx}dy.$$
Notice that from this it follows that
$-i\hbar (d\psi_{qr}/dx)\in L_{2}(\R )$
and the momentum operator is defined at $\psi_{qr}$.
Also (denoting $N_{k}=\sqrt{N}/k$)
$$\theta_{N}=
\sum_{n=0}^{N_{k}-1}(N_{k})^{-1/2}h(n/\sqrt{N})
\phi_{(n/\sqrt{N})((n+1)/\sqrt{N})}$$
converges to $hx\phi_{qr}$ when $\sqrt{N}$ becomes increasingly divisible.
Thus
$$\Phi(\theta_{N})=\sum_{n=0}^{N_{k}-1}(N_{k})^{-1/2}h(n/\sqrt{N})
\psi_{(n/\sqrt{N})((n+1)/\sqrt{N})}$$
converges to $\Phi (hx\phi_{qr})=-i\hbar(d\psi_{qr}/dx)$ (as $\Phi$ is an isometric operator).
Thus $G(-i\hbar (d\psi_{qr}/dx))=\lim_{D}G(\Phi (\theta_{N}))$.
So now it is enough to show that
\begin{equation}\label{eq2}
\lim_{D}G(\Phi (\theta_{N}))=(\xi_{N})/D,
\end{equation}
 where
$$\xi_{N}=\sum_{n=0}^{N_{k}-1}(N_{k})^{-1/2}h(n/\sqrt{N})
v(n),$$
since clearly also $P(G(\psi_{qr}))=(\xi_{N})/D$.

To prove \eqref{eq2}, it is enough to prove it
in the form which we get from \eqref{eq2} by first
applying $G$ to $\Phi (\theta_{N})$
in a direct way and then
replacing
$v(n)$ by
$u(n)$, $n<N$, everywhere.
But then we get a claim, which we have already seen:
Let $f\in L_{2}(\R )$ be such that $f(x)=hxr^{-1/2}$ if $0\le x\le r$
and otherwise $f(x)=0$. Then
letting $\xi^{*}_{N}=\sum_{n=0}^{N_{k}}(N_{k})^{-1/2}h(n/\sqrt{N})
u(n)$,
$F(f)=\lim_{D}\xi^{*}_{N}$
if we define $F$ as it was done first.
On the other hand $F(f)=\lim_{D}F(\theta_{N})$, when we
use the second definition of $F$ 
and of course, $F(\theta_{N})$ is what
we get from $G(\Phi (\theta_{N}))$ if we
replace $v(n)$'s by  $u(n)$'s.
\end{proof}

Next we look at exponentials of the momentum operator.
Let $\a$ and $\b$ be real numbers. Then the subspaces $H^{k}_{N}$
work also for the operators $\a (\b P_{N})^{2}$, $N<\o$.
And the operators $\a (\b P_{N})^{2}$ are self adjoint.
Thus using the
methods from Section~\ref{ultraproducts}, we get an unbounded operator
on $K$.  It is easy to see that this operator is the same
as $\a (\b P)^{2}$. And then still using the  methods from Section~\ref{ultraproducts},
we get the operator $e^{i\a  (\b  P)^{2}}$.

\begin{proposition}
$F\circ e^{i\a (-i\b\hbar (d/dx))^{2}}\circ F^{-1}
=e^{i\a  (\b  P)^{2}}\raj \rng(F)$.
\end{proposition}

\begin{proof}
For all $\psi_{qr}\in S$,
$e^{i\a (-i\b\hbar (d/dx))^{2}}(\psi_{qr})$
can be calculated using the usual series for exponential
and thus by Fact~\ref{fact1.3} \eqref{toinen},
$$(G\circ  e^{i\a (-i\b\hbar (d/dx))^{2}})(\psi_{qr})=
(e^{i\a  (\b  P)^{2}}\circ G)(\psi_{qr}).$$
This shows the claim.
\end{proof}

One can repeat the arguments above (in fact the arguments are easier)
for the operators $x$ and $Q_{N}$ in place of
$-i\hbar (d/dx)$ and $P_{N}$ and get:

\begin{proposition}
$F\circ x\circ F^{-1}=Q\raj F(\dom(x))$ and
$F\circ e^{i\a (\b x)^{2}}\circ F^{-1}
=e^{i\a  (\b  Q)^{2}}\raj \rng(F)$.
\end{proposition}

By the techniques introduced here one can prove similar results
for many other operators than those  studied here
(e.g., for operators $e^{i\a  x}$ and $e^{i\a (-i\hbar (d/dx))}$).

\section{Calculating with eigenvectors}\label{sec:calc}

In this section we show how to use our model and its eigenvectors to calculate the kernel of the time evolution operator. We start by calculating the Feynman propagator for the free particle with the eigenvector method, i.e., as it is often described in physics as the inner product
$$
\langle x_1|K^t|x_0\rangle,
$$
where the $|x_i\rangle$ are eigenvectors of $Q$ and $K^t$ is the time
evolution operator.

The idea behind our method is that we look at eigenvectors
$|x_i\rangle$ in the ultraproduct $H$, find representatives for them
in $H_N$, and do the calculation in our finite dimensional spaces $H_N$.
However, this method leads to a discretising effect stemming from divisibility issues, and the value we get in $H_N$ will be either 0 or $N^{-1/2}\alpha K(x_0,x_1,t)$, where $K(x_0,x_1,t)$ is the propagator (or the kernel of the time evolution operator) as calculated by physicists and $\alpha$ is a real scaling constant (of magnitude $th/m$ for the free particle and $h\sin(\o t)/m\o$ for the harmonic oscillator).

To remedy the discretising effect, one needs to calculate limits of averages over small areas, and this amounts to eventually using the eigenvector method, not to calculate the propagator directly, but to calculate the kernel.

After explaining the kernel calculation for the free particle, we turn to the calculations for the harmonic oscillator, where we need to use a factorisation of the time evolution operator to get nice calculations with eigenvectors.

\subsection{The Feynman propagator of the free particle}\label{propagator}

To calculate the Feynman propagator of the free particle, let us first look at its time evolution operator.

For a time independent
Hamiltonian $H$, the time evolution operator is
$$
K^{t}=e^{-itH/\hbar},
$$
and for a free particle, the Hamiltonian is
$$H=P^{2}/2m,$$
where $m$ is the mass
of the particle (see \cite{Ze};
in \cite{Zi} there is probably a typo here).
Thus the time evolution operator we look at is
$$
K^t=e^{-itP^2/2m\hbar}.
$$

Recall, that the momentum operator $P=-i\hbar\frac{d}{dx}$ is mapped by $F^a_b$ to the operator $P$ defined on $K$, stemming from the operators $P_N$, where
$$
P_N(v(n))=\left\{\begin{array}{ll}
bnN^{-1/2}v(n) & \textrm{when }n<\frac{N}{2}\\
b(n-N)N^{-1/2}v(n)&\textrm{otherwise}\end{array}\right.
$$
In our calculation, however, we will use Gauss sums and thus we need our sums to be of the right form. To attain this, we will instead of $P$ use the operator $P'$ stemming from the operators $P'_N$, where $P'_N(v(n))=bnN^{-1/2}$ for all $n<N$. We can do this, when the operators give the same values for all $v(n)$, i.e., when
$e^{-it(bn)^2/2m\hbar N}=e^{-it(b(n-N))^2/2m\hbar N}$, which is when $tb^2/2mh$ is an integer. 
At first this may look bad but with our choice $b=h$ this 
is not a constant, it comes with a unit  and this unit is the square of
the unit of length (the unit of $h$ is that of
$mass\times  length^{2}/time$). Thus by choosing the unit of length suitably, the requirement can be satisfied.

\begin{lemma}
When $a=1$, $b=h$ and the units are chosen such that $th/2m\in\Z$, then for rational positions $x_0, x_1$ the Feynman propagator in $H_N$ is
$$
\langle x_1\vert K^t \vert x_0\rangle=N^{-1/2}thm^{-1} K(x_0,x_1,t),
$$
when $thm^{-1}$ divides $\sqrt{N}(x_1-x_0)$ and 0 otherwise, 
where $K(x_0,x_1,t)= (m/2\pi i \hbar t)^{1/2}e^{im(x_0-x_1)^2/2\hbar t}$ is the value found by physicists.
\end{lemma}

\begin{proof}
Since we are interested in calculations in almost all $H_N$, we may assume $N$ is large enough such that both $\sqrt{N}x_i$ are integers. To simplify calculations, we will also assume that both $x_0$ and $x_1$ are non-negative. Now (see Notation~\ref{notaatio}) $\ket{x_i}_N=u(\sqrt{N}x_i)$.

Further, as now $tb^2/2mh=th/2m\in\Z$, we can use the operator $P'_N$ that has eigenvalues $hxN^{-1/2}$ for the vectors $v(x)$, and thus as time evolution operator we may use
$$K^t(v(x))=e^{-it(hx)^2/N2m\hbar}v(x)=e^{-i\pi thx^2/mN}v(x), $$
for $0\leq x<N$.

Thus (recall $q=e^{i2\pi /N}$),
$$
\ket{x_0}=\sqrt{\frac{1}{N}}\sum_{n=0}^{N-1}q^{-\sqrt{N}x_0n}v(n)
$$
so 
\begin{align*}
K^t\ket{x_0}&=\sqrt{\frac{1}{N}}\sum_{n=0}^{N-1}q^{-\sqrt{N}x_0n}e^{-i\pi thn^2/mN}v(n)\\
&=\frac{1}{N}\sum_{n=0}^{N-1}\sum_{k=0}^{N-1}q^{-\sqrt{N}x_0n}e^{-i\pi thn^2/mN}q^{nk}u(k),
\end{align*}
and finally
\begin{align*}
\langle x_1|K^t\ket{x_0}
&=\frac{1}{N}\sum_{n=0}^{N-1}q^{-\sqrt{N}x_0n}e^{-i\pi thn^2/mN}q^{n\sqrt{N}x_1}\\
&=\frac{1}{N}\sum_{n=0}^{N-1}e^{-i2\pi\sqrt{N}x_0n/N-i\pi thn^2/mN+i2\pi n\sqrt{N}x_1/N}\\
&=\frac{1}{N}\sum_{n=0}^{N-1}e^{\pi i(2\sqrt{N}(x_1-x_0)n - thm^{-1}n^2) /N}.
\end{align*}

Now we apply number theory:
For integers $c,d,g$, if $cg\ne 0$
and $cg-d$ is even, then
$$\sum_{n=0}^{|g|-1}e^{\pi i(cn^{2}+dn)/g}=|g/c|^{1/2}
e^{\pi i(|cg|-d^{2})/4cg}\sum_{n=0}^{|c|-1}e^{-\pi i(gn^{2}+dn)/c},$$
see \cite{BE} (this result is originally
from Siegel's paper \cite{Si} but at least for us that paper
was hard to find).
So
$$
\langle x_1|K^t\ket{x_0}
=\frac{1}{N}\sqrt{\frac{Nm}{th}}e^{-\pi i/4}e^{\pi i(x_1-x_0)^2/thm^{-1}}\sum_{n=0}^{thm^{-1}-1}e^{-\pi i(Nn^2+2\sqrt{N}(x_1-x_0)n)/thm^{-1}}.
$$
Now the factor outside the sum is (using the fact that $(e^{-i\pi /4})^{2}=1/i$)
$$
N^{-1/2}(m/thi)^{1/2}e^{\pi im(x_1-x_0)^2/th},
$$
which is $N^{-1/2}K(x_0,x_1,t)$. The sum, however, depends on divisibility: If $thm^{-1}$ divides $\sqrt{N}(x_1-x_0)$, the exponent becomes a multiple of $2\pi i$ (we may assume $2thm^{-1}$ divides $N$), so the sum becomes $thm^{-1}$. If $thm^{-1}$ does not divide $\sqrt{N}(x_1-x_0)$, we are summing roots of unity, and the sum amounts to 0.
\end{proof}

At first sight, it seems one could argue, that $\sqrt{N}(x_1-x_0)$ is as divisible as one likes (for large enough $N$). However, if one looks at the ultraproduct, there are, for each position $x$, continuum many orthogonal eigenvectors corresponding to that value. This is because one can choose different sequences of eigenvectors, whose eigenvalues all converge to $x$. And there is no guarantee one has enough divisibility along the way in these sequences. 

I a sense, one wishes to calculate the average of the propagators calculated for all these eigenvectors, and this effect one can get by calculating the kernel instead of the propagator. We show how to do this next.

\subsection{The kernel}\label{kernel}

So let us calculate the kernel
of the time evolution operator  in $H$  and notice that
this way we get the correct answer, i.e., the one calculated by physicists.
We denote by $K(x,y,t)=(m/iht)^{1/2}e^{i\pi m(x-y)^2/ht}$ this correct answer.

Now assuming $\phi$ (and $K(x,y,t)$) are continuous, we have
$$
K(\alpha,\beta,t)=\lim_{\e\to0}\int_{\beta-\e}^{\beta+\e}\int_{\alpha-\e}^{\alpha+\e}\phi(x)K(x,y,t) dx dy\,/\,((2\e)^2\phi(\alpha))
$$
and we wish to calculate this limit in $H$. 

For this, first note that, denoting by $\phi_{[\alpha-\e,\alpha+\e]}(x)$ the function that gets the value $\phi(x)$ when $\alpha-\e\leq x\leq \alpha+\e$, and 0 otherwise,
$$
\int_{\alpha-\e}^{\alpha+\e}\phi(x)K(x,y,t) dx=\int_{\R}\phi_{[\alpha-\e,\alpha+\e]}(x)K(x,y,t) dx=K^t(\phi_{[\alpha-\e,\alpha+\e]})(y).
$$
Second, note that the outer integral in our limit expression can be approximated by Riemann sums, and from this we take the interpretation for 'integrating in $H$', i.e., 
$$
\int_{\beta_1}^{\beta_2}f(y)dy =\lim_{N\to\infty} \sum_{y=\sqrt{N}\beta_1}^{\sqrt{N}\beta_2-1}f(y/\sqrt{N})N^{-1/2}=\lim_D\sum_{y=\sqrt{N}\beta_1}^{\sqrt{N}\beta_2-1}f(y/\sqrt{N})N^{-1/2}.
$$
Now, recalling our embeddings $F_N$ (for the scaling $a=1$, $b=h$):
\begin{align*}
F_{N}(f)= & \sum_{n=0}^{(N/2)-1}N^{-1/4}f(nN^{-1/2})u(n)+\\
 & \sum_{n=N/2}^{N-1}N^{-1/4}f((n-N)N^{-1/2})u(n),
\end{align*}
we have for rational $\beta_i$ (for ease of notation, we assume $0<\beta_1<\beta_2$ and we look at $N$ large enough so that $\sqrt{N}\beta_i$ are integers and $\beta_2<\sqrt{N}/2$)
$$
\sum_{y=\sqrt{N}\beta_1}^{\sqrt{N}\beta_2-1}f(y/\sqrt{N})N^{-1/2}=\sum_{y=\sqrt{N}\beta_1}^{\sqrt{N}\beta_2-1} N^{-1/4}\langle u(y)|F_N(f)\rangle.
$$

Now
$$
F(K^t(\phi_{[\alpha-\e,\alpha+\e]}))=(K_N^t(F_N(\phi_{[\alpha-\e,\alpha+\e]})))_{N<\o}/D
$$
so we can calculate in the models $H_N$, and when $th/2m\in\Z$ (so that we can use $P'$ instead of $P$), and $\alpha$ and $\e$ are rational numbers, we get (as in the propagator case)
\begin{align*}
K_N^t(F_N(\phi_{[\alpha-\e,\alpha+\e]}))
&=K_N^t\left(\sum_{k=\sqrt{N}(\alpha-\e)}^{\sqrt{N}(\alpha+\e)-1}N^{-1/4}\phi(k/\sqrt{N})u(k)\right)\\
&=\sum_{k=\sqrt{N}(\alpha-\e)}^{\sqrt{N}(\alpha+\e)-1}N^{-1/4}\phi(k/\sqrt{N})N^{-1}\sum_{n=0}^{N-1}\sum_{l=0}^{N-1}e^{2\pi il(n-k)/N}e^{-i\pi thl^2/mN}u(n).
\end{align*}
Thus to calculate the integral
$$
\int_{\beta-\e}^{\beta+\e}K^t(\phi_{[\alpha-\e,\alpha+\e]})(y)dy
$$
we can calculate in $H_N$, assuming $th/2m$ is an integer and $\alpha,\beta,\e$ are rational numbers:
\begin{gather*}
\sum_{n=\sqrt{N}(\beta-\e)}^{\sqrt{N}(\beta+\e)-1}N^{-1/4}\langle u(n)|K^t_N(F_N(\phi_{[\alpha-\e,\alpha+\e]}))\rangle \\
=\sum_{n=\sqrt{N}(\beta-\e)}^{\sqrt{N}(\beta+\e)-1}N^{-1/4} \sum_{k=\sqrt{N}(\alpha-\e)}^{\sqrt{N}(\alpha+\e)-1}N^{-1/4}\phi(k/\sqrt{N})N^{-1}\sum_{l=0}^{N-1}e^{2\pi il(n-k)/N}e^{-i\pi thl^2/mN}\\
=\sum_{n=\sqrt{N}(\beta-\e)}^{\sqrt{N}(\beta+\e)-1}\sum_{k=\sqrt{N}(\alpha-\e)}^{\sqrt{N}(\alpha+\e)-1}N^{-3/2}\phi(k/\sqrt{N})\sum_{l=0}^{N-1}e^{\pi i (2(n-k)l- thm^{-1}l^2)/N}.
\end{gather*}

Now we may use the Gaussian sum formula as before. We have two
cases: If $(n-k)$ is divisible by $th/m$, the last sum becomes
$$
\sqrt{N}\sqrt{\frac{m}{thi}}e^{\pi i (n-k)^2/Nthm^{-1}}\cdot thm^{-1}=\sqrt{N}thm^{-1}K(n/\sqrt{N},k/\sqrt{N},t).
$$
However, if $thm^{-1}$ does not divide $(n-k)$, we end up summing
roots of unity, and the result is 0. So all in all we are summing over $(\sqrt{N}2\e)^2$ pairs of coordinates, and on average one in $thm^{-1}$ many gives a value that is scaled by a factor of $thm^{-1}$. So on average the sum we get is
$$
\sum_{n=\sqrt{N}(\beta-\e)}^{\sqrt{N}(\beta+\e)-1}
\sum_{k=\sqrt{N}(\alpha-\e)}^{\sqrt{N}(\alpha+\e)-1}
N^{-1} \phi(k/\sqrt{N}) K(k/\sqrt{N},n/\sqrt{N},t).
$$
Now what we wanted to calculate was
$$
\lim_{\e\to 0}\lim_D
\sum_{n=\sqrt{N}(\beta-\e)}^{\sqrt{N}(\beta+\e)-1}
\sum_{k=\sqrt{N}(\alpha-\e)}^{\sqrt{N}(\alpha+\e)-1}
N^{-1} \phi(k/\sqrt{N}) K(k/\sqrt{N},n/\sqrt{N},t)\,/\,(2\e)^2\phi(\alpha)
$$
and for a small enough $\e$ this is approximately
$$
\lim_D(\sqrt{N}2\e)^2 N^{-1}\phi(\alpha)K(\alpha,\beta,t)/(2\e)^2\phi(\alpha)=K(\alpha,\beta,t).
$$

Thus we see that the time evolution operator gathers weight into
'spikes', and calculating the propagator as an inner product picks out
these spikes, whereas one gets the correct value by averaging these
spikes over a small area around the spike.

\subsection{The kernel of the harmonic oscillator}\label{oscillator}

Next we turn our attention to the harmonic oscillator. 
The Hamiltonian for
the system is
$$
H=\frac{P^2}{2m}+\frac{1}{2}m\omega^2Q^2,
$$
where $\omega$ is the angular frequency.
Thus the time evolution operator one gets from the Schr\"odinger equation
for a time independent Hamiltonian is
$$
K^t=e^{-itH/\hbar}=e^{-\frac{it}{\hbar}(\frac{P^2}{2m}+\frac{1}{2}m\omega^2Q^2)}.
$$

To be able to use our eigenvector method to calculate the Feynman propagator $\langle x_1|K^t|x_0\rangle$ or the kernel $K(x_0,x_1,t)$ using the bases in the space $H_N$, we need
to know the matrix of the time evolution operator in one
of our bases. To find such a matrix is hard. One solution
used by the physicists,
is that they factorise the operator. In particular, they try to find
a factorisation that separates the $P$- and $Q$-parts of the
time evolution operator. The literature knows various factorisations of the
operator, stated for operators satisfying the commutator relation
$[Q,P]=i\hbar$. 
A first obstacle for us is the already mentioned fact that no finite
dimensional Hilbert space can satisfy this relation. However, we reason as follows: As we can embed the standard $L_2$-space into $H$ and the factorisation holds in $L_2$, it holds in the image of the embedding. More specifically, what we calculate is not the propagator in the spaces $H_N$, and the limit may not be the propagator in the ultraproduct model $H$ either (we cannot guarantee the right commutator $[Q,P]$ outside the image of $L_2(\R)$, as - although the Weyl commutator relation holds in a dense set - the Weyl operators do not form continuous representations in all of $H$. However, what we are calculating, is the kernel, and thus we can use any expression for $K^t_N$ as long as it satisfies
$$
F(K^t(\phi))=(K_N^t(F_N(\phi)))_{N<\o}/D
$$
for a dense set of functions $\phi$, and this will be guaranteed by our observations in section~\ref{embedding}.

\medskip

We look at a factorisation of the time evolution operator from \cite{QA}. This
is formulated for operators $A$, $B$ and $C$
satisfying $[A,B]=C$, $[A,C]=2\gamma A$, $[B,C]=-2\gamma B$ (and these
are true of the operators we are interested in, when $[Q,P]=i\hbar$
and Baker-Campbell-Hausdorff holds for the operators),
and it gives coefficients $\alpha,\beta$ satisfying
$$
e^{A+B}=e^{\alpha A}e^{\beta B}e^{\alpha A}.
$$
Quijas and Aguilar seem to have some typos in their paper (the third
term of the third equation in (55) should be positive, and (59) is not
what you get by substitution). But their result satisfies the equations
given by the method, producing 
$$
\alpha=\frac{1}{\sqrt{\gamma}}\tan(\sqrt{\gamma}/2)
$$
and
$$
\beta=\frac{1}{\sqrt{\gamma}}\sin(\sqrt{\gamma}).
$$

Now if we substitute $A=-\frac{itm\omega^2}{2\hbar}Q^2$, $B=-\frac{it}{2m\hbar}P^2$
and $C=-\frac{it^2\omega^2}{2\hbar}(QP+PQ)$, this satisfies the
commutator requirements with $\gamma=t^2\omega^2$ and gives
$$\alpha=\frac{1}{t\omega}\tan(t\omega/2)$$
and 
$$\beta=\frac{1}{t\omega}\sin(t\omega)$$
and thus
\begin{align*}
K^t &=e^{-itH/\hbar}=e^{-\frac{it}{\hbar}(\frac{P^2}{2m}+\frac{1}{2}m\omega^2Q^2)}\\
&=e^{-\frac{im\omega\tan(t\omega/2)}{2\hbar}Q^2}e^{-\frac{i\sin(t\omega)}{2\omega m\hbar}P^2}e^{-\frac{im\omega\tan(t\omega/2)}{2\hbar}Q^2}.
\end{align*}

Now we can write this as
$$
K^t=e^{irQ^2}e^{isP^2}e^{irQ^2},
$$
where $r=-\frac{m\omega\tan(t\omega/2)}{2\hbar}$ 
and $s=-\frac{\sin(t\omega)}{2\omega m\hbar}$, 
and recall that $e^{irQ^2}$ corresponds to the operator $u(x)\mapsto e^{ir\underline{x}^2/N}u(x)$ and $e^{isP^2}$ to the operator $v(x)\mapsto e^{is(h\underline{x})^2/N}v(x)$, where
$$
\underline{x}=\left\{\begin{array}{ll}
x &\textrm{if }x<\frac{N}{2},\\
x-N&\textrm{if }\frac{N}{2}\leq x<N
\end{array}\right.
$$
Thus we can calculate, by switching between the bases $\{u(x)\vert\ x<N\}$ and
$\{v(x)\vert\ x<N\}$,
$$
K_N^t(u(x))=e^{ir\underline{x}^2/N}\frac{1}{N}\sum_{n=0}^{N-1}\sum_{k=0}^{N-1}q^{n(k-x)}e^{ish^2\underline{n}^2/N}e^{ir\underline{k}^2/N}u(k).
$$

Now we calculate the kernel as for the free particle, approximating the integral
$$
\int_{\beta-\e}^{\beta+\e}\int_{\alpha-\e}^{\alpha+\e}\phi(x)K(x,y,t) dx dy=\int_{\beta-\e}^{\beta+\e}K^t(\phi_{[\alpha-\e,\alpha+\e]})dy,
$$
where, as before, $\phi_{[\alpha-\e,\alpha+\e]}$ denotes the restriction of $\phi$ to the interval $[\a-\e,\a+\e]$. So we calculate (for large enough $N$ and positive, rational $\a$, $\b$ and $\e$)
\begin{equation}\label{apprint}
\begin{split}
\sum_{k=\sqrt{N}(\b-\e)}^{\sqrt{N}(\b+\e)-1}N^{-1/4}\left\langle u(k)\left|\, \sum_{l=\sqrt{N}(\a-\e)}^{\sqrt{N}(\a+\e)-1}N^{-1/4}\phi(lN^{-1/2})K^t_N(u(l))\right\rangle\right.\\
=
\sum_{k=\sqrt{N}(\b-\e)}^{\sqrt{N}(\b+\e)-1} N^{-1/4}\sum_{l=\sqrt{N}(\a-\e)}^{\sqrt{N}(\a+\e)-1}N^{-1/4}\phi(lN^{-1/2})e^{ir\underline{l}^2/N}\frac{1}{N}\sum_{n=0}^{N-1}q^{n(k-l)}e^{ish^2\underline{n}^2/N}e^{ir\underline{k}^2/N}\\
=
\sum_{k=\sqrt{N}(\b-\e)}^{\sqrt{N}(\b+\e)-1} \sum_{l=\sqrt{N}(\a-\e)}^{\sqrt{N}(\a+\e)-1}N^{-3/2}\phi(lN^{-1/2}) e^{ir(l^2+k^2)/N}\sum_{n=0}^{N-1}q^{n(k-l)}e^{ish^2\underline{n}^2/N}.
\end{split}
\end{equation}
Now the last sum is of the form
$$
\sum_{n=0}^{N-1}q^{n(k-l)}e^{ish^2\underline{n}^2/N}
=\sum_{n=0}^{N-1}e^{\pi i (2n(k-l)-h\sin(t\omega)\underline{n}^2(\omega m)^{-1})/N},
$$
and when (we have scaled units such that) $\frac{h\sin(\omega t)}{\omega m}$ is an integer, we can switch $\underline{n}$ to $n$. Thus we can use the same Gauss formula as before, with $g=N$, $c=-h\sin(\omega t)/(\omega m)$ and $d=2(k-l)$, and get
\begin{gather*}
\sum_{n=0}^{N-1}e^{\pi i (2n(k-l)-h\sin(t\omega)\underline{n}^2(\omega m)^{-1})/N}\\
=
N^{1/2}\left(\frac{\omega m}{h\sin(\omega t)}\right)^{1/2}e^{\pi i(-\frac{1}{4}+(k-l)^2\omega m (h\sin(t\omega))^{-1}/N)}\sum_{n=0}^{h\sin(\omega t)/(\omega m)-1}e^{\pi i (Nn^2+2(k-l)n)/(h\sin(\omega t)/(\omega m))}\\
= N^{1/2}\left(\frac{\omega m}{ih\sin(\omega t)}\right)^{1/2}e^{\frac{\pi i(\frac{k}{\sqrt{N}}-\frac{l}{\sqrt{N}})^2\omega m}{h\sin(t\omega)}}\sum_{n=0}^{h\sin(\omega t)/(\omega m)-1}e^{\pi i (Nn^2+2(k-l)n)\omega m/(h\sin(\omega t))}.
\end{gather*}
And similarly to before, the last sum is $\frac{h\sin(\omega t)}{\omega m}$ when $k-l$ is divisible by $\frac{h\sin(\omega t)}{\omega m}$, and 0 otherwise, and the divisible case occurs for one in $\frac{h\sin(\omega t)}{\omega m}$ of the pairs $(k,l)$. If we calculate the average value we thus get 
\begin{gather*}
\sum_{k=\sqrt{N}(\b-\e)}^{\sqrt{N}(\b+\e)-1} \sum_{l=\sqrt{N}(\a-\e)}^{\sqrt{N}(\a+\e)-1}N^{-3/2}\phi(lN^{-1/2}) e^{ir(l^2+k^2)/N}\sum_{n=0}^{N-1}q^{n(k-l)}e^{ish^2\underline{n}^2/N}
\\
=
\sum_{k=\sqrt{N}(\b-\e)}^{\sqrt{N}(\b+\e)-1} \sum_{l=\sqrt{N}(\a-\e)}^{\sqrt{N}(\a+\e)-1}N^{-3/2}\phi(\frac{l}{\sqrt{N}}) e^{-\frac{\pi im\omega\tan(t\omega/2)((\frac{l}{\sqrt{N}})^2+(\frac{k}{\sqrt{N}})^2)}{h}}
N^{1/2}\\
{} \left(\frac{\omega m}{ih\sin(\omega t)}\right)^{1/2}
e^{\frac{\pi i(\frac{k}{\sqrt{N}}-\frac{l}{\sqrt{N}})^2\omega m}{h\sin(t\omega)}}
\\
=
\sum_{k=\sqrt{N}(\b-\e)}^{\sqrt{N}(\b+\e)-1} \sum_{l=\sqrt{N}(\a-\e)}^{\sqrt{N}(\a+\e)-1} N^{-1} \phi(\frac{l}{\sqrt{N}}) 
K_{\textrm{phys}}\left(\frac{l}{\sqrt{N}},\frac{k}{\sqrt{N}},t\right),
\end{gather*}
where $K_{\textrm{phys}}(x,y,t)=\left(\frac{m\omega}{ih\sin(\omega t)}\right)^{1/2}e^{-\frac{\pi m\omega(cos(t\omega)(x^2+y^2)-2xy)}{ih\sin(\omega t)}}$ is the value physicists get for the kernel (or propagator). In the last equation we have used the fact that $\tan(x/2)=(1-\cos(x))/\sin(x)$.

Again, for a small enough $\e$, this yields approximately
$$
(\sqrt{N}2\e)^2 N^{-1} \phi(\alpha)K_{\textrm{phys}}(\alpha,\beta,t)=4\e^2\phi(\alpha)K_{\textrm{phys}}(\alpha,\beta,t).
$$
Thus the kernel we calculate is 
$$
K(\alpha,\beta,t)=\lim_{\e\to0}\frac{1}{(2\e)^2\phi(\alpha)}\lim_D 4\e^2\phi(\alpha)K_{\textrm{phys}}(\alpha,\beta,t)=K_{\textrm{phys}}(\alpha,\beta,t).
$$

\section{Other scalings}\label{scalings}

Both the 'renormalisation' factor $N^{-1/2}$ and the discretisation effect one gets when mapping a continuous function to a finite-dimensional vector come up in all approaches with finite-dimensional approximations. However, we have not encountered the additional discretisation stemming from the Gauss formula in other works. The reason seems to be the tradition of scaling units heavier than we do, e.g., \cite{TC} use the same Gauss formula, but there units are scaled so that the coefficient giving the period ($thm^{-1}$ in the case of the free particle) is 1. But such a scaling means one is using a different model for each studied time difference. Although our method requires some scaling to be able to calculate with integer values, we can find a minimum grid (both in time and space), so that any multiple time slots and multiple position differences of these minimum distances, can be calculated in the same model.

Another approach, that at first seems to remove the divisibility issue, is scaling the operators. Instead of choosing $a=1$ and $b=h$, as was done in this paper, one can for the free particle choose $a=thm^{-1}$. Then taking this into account both for the operators used and for the eigenvectors, one calculates essentially as in this paper (scaling units such that $a=thm^{-1}$ is the inverse of an integer to be able to switch between $P$ and $P'$ and use the Gauss sum formula), and gets for rational $x_0, x_1$ in $H_N$ for large enough $N$,
$$
\langle x_{1}|K^{t}|x_{0}\rangle =
N^{-1/2}(m/thi)^{1/2}e^{\pi i m (x_1-x_0)^2/th}
=N^{-1/2}K_{\textrm{phys}}(x_0,x_1,t),
$$
i.e., $N^{-1/2}$ times the value from physics literature. However, although this value looks right (up to a removal of the 'renormalisation' factor $N^{-1/2}$), it is wrong from the point of view of the model, as it gives the wrong probabilities:
If $\phi$ is the state
(in $L_2$) $a^{-1/2}N^{1/4}\chi_{[x_0,x_0+\frac{a}{\sqrt{N}}]}$ (where $\chi$
denotes the characteristic function), then the scaled version of $F_N$
maps it onto $|x_0\rangle$. If we calculate the probability of a
particle starting in state $\phi$ at time 0 to end up in
$[x_1,x_1+\frac{a}{\sqrt{N}}]$ in time $t$, in $H_N$ this is given by
$$
|\langle x_1|K^t|x_0\rangle|^2=N^{-1}K_{\textrm{phys}}(x_0,x_1,t)^2.
$$
However, in the standard model the corresponding probability is
\begin{align*}
&\int_{[x_1,x_1+\frac{a}{\sqrt{N}}]}|\int_{\R}\phi K(x,y,t)dx|^2 dy \\
&=\int_{[x_1,x_1+\frac{a}{\sqrt{N}]}}|\int_{[x_0,x_0+a/\sqrt{N}]}a^{-1/2}N^{1/4}
K(x,y,t)dx|^2 dy 
\end{align*}
and this is approximately
$$
\frac{a}{\sqrt{N}}|\frac{a}{\sqrt{N}}a^{-1/2}N^{1/4}K(x_0,x_1,t)|^2=a^2N^{-1}K(x_0,x_1,t)^2,
$$
so with the scaling $a=thm^{1}$  we get different probabilities in $H_N$ and $L_2(\R)$.

Note that the propagator we calculated in \ref{propagator} gives the right probabilities, as long as one does not use it for point values: if we fix two intervals and calculate the probability of a particle in one to end up in the other in a given time, the value we get in $H_N$ (when transferring the situation from $L_2(\R)$ to $H_N$ via $F_N$) will, when $N$ grows, approach the value calculated in the standard model. Thus the propagator is correct, although it does not look so at first sight: both the discretisation factor and 'renormalisation' factor are necessary for correct probabilities.

Finally, although scaling doesn't help in getting rid of the discretisation, it opens up possibilities for applications of our model. We conclude with an example from chemistry.

\subsection{On Bohr's model for hydrogen atom}\label{bohr}

In this section we look at what happens in our model if we
radically change the scaling $a'$ and $b'=Nh/a'$ from section~\ref{building}.
This time we let $a'=a$ be a positive real and then
$b'=Nh/a=bN$. Then in $H_{N}$ we have eigenvectors $u(n)$
for positions $an/N$ if $n<N/2$
and $a(n-N)/N$ otherwise. So when $N$ goes to infinity,
the model covers the interval $[-a/2,a/2[$, and there is a natural embedding of $L_2([-a/2,a/2])$ into our ultraproduct. Defining $Q$ and $P$ in the standard way in $L_2([-a/2,a/2])$ gives the commutation relation $[Q,P]=i\hbar$, whenever it is defined. However, to make $P$ self adjoint in this space we need to restrict it to the space of continuously differentiable functions $f$ on $[-a/2,a/2]$ with the boundary condition $f(-a/2)=f(a/2)$. Defining $e^{itP}$ then gives an operator consisting of 'translation with wraparound' (see Example 14.5 of \cite{Ha}). Thus the natural interpretation of our model is not $L_2([-a/2,a/2])$, but the $L_2$-space over a circle.

 So, in our model, we interpret the positions as a circle around the origin in $\R^3$.
More precisely, for a fixed $a$,
$u(n)$ is an eigenvector for the position
$$
cc(n)=(r\ \cos(2\pi \underline{n}/N), r\ \sin(2\pi \underline{n}/N),0)
$$
where
$r=a/2\pi$ is the radius of the circle, and, as before, $\underline{n}$ denotes $n$, when $n<N/2$ and $n-N$ otherwise.

To be able to consider circles of various radii, we think of the model as having scaled operators $P_a$ and $Q_a$ for all $a$ (and $b$ such that $ab=h$). Recall that with this scaling $Q_a(u(n))=a\underline{n}/N\ u(n)$ and $P_a(v(n)=b\underline{n}v(n)$.
In order to simplify the notation,  we write $v(\underline{n})$ for $v(n)$, i.e., talk about vectors $v(k)$ for $k=-N/2,\dots,N/k-1$,  and note that the sign tells the direction of the momentum (vector). Also we denote by 
$S_{r}$ the circle $\{ (r\ cos(x),r\ sin(x),0)\vert\ x\in [0,2\pi [\}$.

We aim to look at Bohr's model for the hydrogen atom,
and for this we consider the proton being in the origin
and $H_N$ describing the state of the electron.
So suppose that the state
$\psi$ is $v(n)$. Then the momentum of the electron is $bn=n\hbar /r$
(keep in mind that $ab=h$ and $a=2\pi r$). Thus, proceeding as Bohr, we can calculate its kinetic energy
$$
E_{kinetic}=\frac{p^2}{2m}=\frac{n^2\hbar^2}{2m_er^2},
$$
where $m_e$ is the  mass of the electron,
and its Coulomb potential
$$
E_{potential}=-\frac{e^2}{4\pi\varepsilon_0 r},
$$
where $e$ is the elementary charge, and $\varepsilon_0$ is the permittivity of vacuum.
In calculations like these, one usually uses atomic units, where the numerical values of $m_e$, $e$, $\hbar$ and $1/(4\pi\varepsilon_0)$ are all 1. Thus we obtain the total energy (in atomic units)
$$
E_{kinetic}+E_{potential}=\frac{n^2}{2r^2}-\frac{1}{r}.
$$
And then we can determine
$a=2\pi r$, by minimising the total energy
$n^{2}/2r^{2}-1/r$. This is exactly what Bohr did
and an easy calculation gives $r=n^{2}$ and thus the energy
is $-1/2n^{2}$, which is in agreement with the Bohr model.

\begin{remark}
In the literature one usually does not minimise the energy
but instead  finds $r$ by requiring that the Coulomb force
and the 'centrifugal force' cancel each other out.
Now the Coulomb force is $-e^2/(4\pi\varepsilon_0 r^2)$, i.e., $-1/r^2$ in atomic units, and the 'centrifugal force' is $p^2/m_er$, which in our model is $(bn)^2/m_er=\hbar^2n^2/m_er^3)$, i.e., $n^2/r^3$ in atomic units.
But notice that $-(n^{2}/r^{3}-1/r^{2})$ is nothing but the derivative
of the energy $n^{2}/2r^{2}-1/r$ and thus the two methods give the
same results.
\end{remark}

One can ask, why we only looked at states of the form $v(n)$ above.
Of course states like $(1/\sqrt{2})(v(n)+v(-n))$ make perfect sense but in the hydrogen atom case
these states give the same results as $v(n)$.
However states
like $(1/\sqrt{2})(v(2)+v(1))$ do not make sense
in our model, because
$v(2)$ and $v(1)$ give different values for $a$ (different radii), and in our model we can have only
one $a$ at a time.

\section*{Acknowledgements}
The authors wish to thank J.~Lukkarinen for his help with kernels of time evolution operators.

\end{document}